\documentclass{amsart}
\usepackage{amssymb,latexsym,amsmath}
\newtheorem{theorem}{Theorem}[section]
\newtheorem{lemma}{Lemma}[section]
\newtheorem{proposition}{Proposition}[section]
\newtheorem{definition}{Definition}[section]

\newtheorem{remark}{Remark}[section]
\numberwithin{equation}{section}

\begin{document}
\title{Conformal Graph Directed Markov Systems }
\author{Andrei E. Ghenciu}
\date{02/20/2007}
\address{Department of Mathematics and Statistics\\
University of Alaska at Fairbanks\\
Fairbanks, Alaska, 99775, U.S.A.}
\email{ffeag@uaf.edu}
\author{Dan Mauldin}
\date{11/01/2007}
\address{Department of Mathematics\\
University of North Texas\\
Denton,Texas, 76203, U.S.A.}
\email{mauldin@unt.edu}
     
\begin{abstract}
We present the main concepts and results for Graph Directed Markov Systems that have a finitely irreducible incidence matrix. We then see how these results change when the incidence matrix is not assumed to be finitely irreducible.
\end{abstract} 
\maketitle 
\section{Introduction}
The theory of the limit set generated by the iteration of finitely many similarity maps
has been well developed for some time now. A more complicated theory of the limit set generated 
by the iteration of infinitely many uniformly contracting conformal maps was developed by Mauldin and Urbanski 
in \cite{MU2}.
Several years after that, they explored the geometric and dynamic properties of a far reaching generalization of conformal iterated function systems, called Graph Directed Markov Systems (GDMS's) 
(see \cite{MU4}).
\\
\\
Several concepts are at the core of the analysis of a conformal Graph Directed Markov System (CGDMS); among them: the topological pressure function, 
the Hausdorff dimension of the limit set and the conformal measure supported on the limit set. The connections between these concepts have been intensively studied by Mauldin and Urbanski, specially in the case when the incidence matrix of the CGDMS is finitely irreducible (see \cite{MU4}). In this paper we focus our attention on the most general case when the incidence matrix is not necessarily finitely irreducible.
\\ 
\\
The paper is organized as follows. In the second section we present the main concepts, constructions and results for CGDMS's as they are presented in \cite{MU4}.
We show how one can construct the limit set by performing an infinite directed walk through the graph. This leads to a natural map from the coding space to the points of the limit set. We look at several properties that the incidence matrix can have. The most important case is when the incidence matrix is finitely irreducible; in this situation most of the results from the theory of the CIFS's ca be carried over. We end this section with the definition of the conformal measure supported on the limit set of a CGDMS.\\
\\
The third section includes the most important results for a CGDMS with its incidence matrix finitely irreducible. In the fourth section we analyze how these results change in the most general case when the incidence matrix is not supposed to be finitely irreducible anymore. Even if we start with 
a CGDMS whose incidence matrix is finitely irreducible, one likes to study the subsystems of the original system. The incidence matrix restricted to the subsystem may not be finitely irreducible.

\newpage
\section{Preliminaries}
To introduce the Graph Directed
Markov Systems we need a directed multigraph $(V,E,i,t)$ and an associated
incidence matrix $A$. The multigraph consists of a finite
set $V$ of vertices and a countable (possibly infinite) set of directed edges $E$ and
two functions $i,t:E \to V$. For each edge $e$, $i(e)$ is
the initial vertex of the edge and $t(e)$ is the terminal vertex of
$e$. We also have a function $A:E^2 \to
\{0,1\}$ called an incidence matrix. The matrix $A$ tells us which edge may follow a given edge. We
define the set of admissible words by:
\[
E^{\infty}_A=\{\omega = (\omega_1,\omega_2,...,\omega_i,...) \in E^{\infty} : A_{\omega_i\omega_{i+1}}=1\ \forall i\geq 1\}.
\]
The set of all finite subwords of $E^{\infty}_A$ will be denoted by $E^*_A$.
By $E^n_A$ we denote the set of all subwords of $E^{\infty}_A$ of length $n\geq1$. If $\omega \in E_A^*$,
by $|\omega|$ we mean the length of the word $\omega$.
If $\omega \in E^{\infty}$ and $n \geq1$, then
\[
\omega|_n=\omega_1\omega_2...\omega_n.
\]
A Graph Directed Markov System (GDMS) consists of a directed multigraph, a set of non-empty compact metric spaces $\{X_v \}_{v\in V}$, a number s,
$0<s<1$, and for every
$e \in E$, a 1-to-1 contraction
$\varphi_e : X_{t(e)} \to X_{i(e)}$ with
a Lipschitz constant $\leq s$ and an incidence matrix $A$. In short, the set
\[
S=\{\varphi_e:X_{t(e)} \to X_{i(e)}\}_{e\in E}
\]
is called a GDMS. 
\\
\\
For each $\omega=(\omega_1,\omega_2,...,\omega_n)$ consider the map coded by $\omega$:
\[
\varphi_{\omega}=\varphi_{\omega_1} \circ ... \circ \varphi_{\omega_n}:X_{t(\omega_n)} \to X_{i(\omega_1)}.
\]
For $\omega \in E^{\infty}_A$, the sets $\{\varphi_{\omega_n}(X_{t(\omega_n)})\}_{n\geq1}$ form a
descending sequence of non-empty compact subsets of $\varphi_{\omega_1}(X_{t(\omega_1)})$.
\\
\\
Since for every $n\geq1$,
$diam(\varphi_{\omega_n}(X_{t(\omega_n)})
\leq s^n diam(X_{t(\omega_n)})
\leq s^n max\{diam(X_v):v\in V \}$, we conclude that the intersection:
\[
\bigcap_{n \geq 1}\varphi_{\omega_n}(X_{t(\omega_n)})
\]
is a singleton and we denote its only element by $\pi(\omega)$.
Thus we can define the map:
\[
\pi:E^{\infty}_A \to X
\]
from $E^{\infty}_A$ to $X$, where $X$ is the disjoint union of the compact sets $X_v$. We will call this map
the \emph{coding map}. The set
\[
J=J_{E,A}=\pi(E^{\infty}_A)
\]
will be called the \emph{limit set} of the GDMS $S$.
From now on we assume that $\forall a \in E$ there exists $b \in E$ so that $A_{ab}=1$.
Otherwise, if there exists $a \in E$ so that $A_{ab}=0$ for every $b \in E$ then the limit set
$J_{E,A}$ would be the same as the limit set
$J_{E \setminus \{b\},A}$
(where $A$ is defined now on $(E \setminus \{b\})^2$).
\\
\\
If the set of vertices of a GDMS is a singleton and all the entries in the associated incidence
matrix are 1, then the GDMS is called an \emph{iterated function system}.
\\
\\
\newpage
We want to emphasize that we have two directed graphs that play an important role in our study. The first one is the given graph $(V,E,i,t)$. The second one,
$G_{E,A}$, which we will call the edge directed graph, is determined by $E$. The vertices of the second directed graph are the edges of the original one. The second graph may have infinitely many vertices and edges.\\
\\
The incidence matrix $A$ is said to be \emph{irreducible} if for any two edges $i,j \in E$ there exists a path
$\omega \in E_A^*$ so that $\omega_1=i$ and $\omega_{|\omega|}=j$. This is equivalent to saying that the
directed graph $G_{E,A}$ is \emph{strongly connected}: for any two vertices there exists a path starting at one and
ending at the other.\\
\\
We say that $A$ is \emph{primitive} if there exists $p \geq1$ such that all the entries of $A^p$ are positive
or, in other words, for any two edges $i,j \in E$ there exists a path
$\omega \in E_A^p$ so that $\omega_1=i$ and $\omega_{|\omega|}=j$.\\
\\
The matrix $A$ is said to be \emph{finitely irreducible} if there exists a finite set $H \subset E_A^*$
so that for any two edges $i,j \in E$ there exists a path $\omega \in H$ so that $i\omega j \in E_A^*$.
\\
\\
The matrix $A$ is said to be \emph{finitely primitive} if there exists a finite set $H \subset E_A^*$
consisting of words of the same length such that for any two edges $i,j \in E$ there exists a path $\omega \in H$ so that $i\omega j \in E_A^*$.
\\
\\
We call a GDMS \emph{conformal} (CGDMS) if the following conditions are satisfied:\\
\\
(1)\ \ For every $v\in V$, $X_v$ is a compact connected subset of a Euclidean space $\mathbf{R}^d$ (the
dimension d common for all vertices) and $X_v=Cl(Int(X_v))$.\\
\\
(2)\ \ (Open Set Condition)(OSC)\\
For every a,b in E, $a\ne b$
$\varphi_a (Int(X_{t(a)})\bigcap\varphi_b (Int(X_{t(b)}) =\O.$\\
\\
(3)\ \ For every vertex $v \in V$ there exists an open connected set $W_v$, $X_v \subset W_v$
so that for every $e \in I$ with $t(e)=v$, the map $\varphi_e$ extends to a $C^1$ conformal
diffeomorphism of $W_v$ into $W_{i(e)}$.\\
\\
(4)\ \ (Cone property) There exists $\gamma,l >0$, $\gamma < \frac{\pi}{2}$ such that for every
$x \in X$ there exists an open cone $Con(x,\gamma,l) \subset Int(X)$ with vertex x, central
angle of measure $\gamma$, and altitude $l$.\\
\\
(5)\ \ There are two constants $L \geq 1$ and $\alpha > 0$ so that:
\[
||\varphi_e'(y)| - |\varphi_e'(x)|| \leq L ||( \varphi_e')^{-1}||^{-1} ||y - x||^{\alpha}
\]
for every $e \in E$ and for every pair of points $x,y \in X_{t(e)}$, where
$|\varphi_e'(x)|$ represents the norm of the derivative. This says that the norms of the derivative
maps are all Holder of order $\alpha$ with the Holder constant depending on the maps.\\

\begin{remark}
If $d \geq 2$ and a family $S=\{\varphi_e:X_{t(e)} \to X_{i(e)}\}_{e\in E}$ satisfies conditions
(1) and (3) from the definition of a GDMS being conformal then it also satisfies (5)
for the same definition with $\alpha=1$.
\end{remark}
As a straightforward consequence of (5) we get the following:\\
\\
(6)\ \ (Bounded Distortion Property) There exists $K \geq 1$ such that for all $\omega \in E_A^*$
and for all $x,y \in X_{t(\omega)}$
\[
|\varphi_{\omega}'(y)| \leq K |\varphi_{\omega}'(x)|.
\]
\\
Given
$t \geq 0$,
a probability measure $m$
is said to be $t-conformal$ provided it is supported
on the limit set $J$ and the following two conditions are satisfied:\\
\\
For every $\omega \in E_A^*$
and for every Borel set $F \subset X_{t(\omega)}$
\[
m(\varphi_{\omega}(F)) = \int_F{|\varphi_{\omega}'|}^t dm
\]
and for all incomparable words $\omega,\tau \in E_A^*$
\[
m(\varphi_{\omega }(X_{t(\omega )}) \bigcap \varphi_{\tau }(X_{t(\tau )})) = 0
\]
The first condition is a change of variable formula and the second is
a measure theoretic open set condition.\\
A simple inductive argument shows that it is enough to require these conditions to
hold at the first level:\\
For every edge $e \in E$ and
for every Borel set $F \subset X_{t(e)}$
\[
m(\varphi_e (F)) = \int_F{|\varphi_e '|}^t dm
\]
and for all edges $a,b \in E$, $a \ne b$
\[
m(\varphi_{t(a)}(X_{t(a)}) \bigcap \varphi_{t(b)}(X_{t(b)})) = 0.
\]
\newpage
\section{The Topological Pressure Function and Finitely Irreducible CGDMS's}
We now define the topological pressure function, a central object in the theory of Conformal Graph Directed Markov Systems.\\
\\
Given $t \geq 0$ and $n \geq 1$ we denote:
\[
Z_{n,E,A}(t) = \sum_{\omega \in E_A^n} ||\varphi_{\omega}'||^t.
\]
Let $\theta_{E,A} = inf \{t \geq 0 \ \ | \ \ Z_{n,E,A}(t) < \infty$ for all but finitely many n's $\}$. 
\\
If $t > \theta_{E,A}$ we put: 
\[
P_{E,A}(t) = \lim_{n \to \infty} \frac{1}{n}\ \ \ln Z_{n,E,A}(t) = \inf_{n \to \infty} \frac{1}{n}\ \ \ln Z_{n,E,A}(t).
\]
If $t < \theta_{E,A}$, we define $P_{E,A}(t) = \infty$.\\
If $Z_{n,E,A}(\theta_{E,A}) < \infty$ for all but finitely many n's we put:
\[
P_{E,A}(\theta_{E,A}) = \lim_{n \to \infty} \frac{1}{n}\ \ \ln Z_{n,E,A}(\theta_{E,A}) = \inf_{n \to \infty} \frac{1}{n}\ \ \ln Z_{n,E,A}(\theta_{E,A}).
\]
Otherwise, $P_{E,A}(\theta_{E,A}) = \infty$.\\ 
We call $P_{E,A}(t)$ the \emph{topological pressure} function.\\
\\
For every $n \geq1$ let $\theta_{n,E,A}=\inf \{t \ \ | \ \ Z_{n,E,A}(t) < \infty \}$
\\
\\
The next proposition is a slightly different version of Proposition 4.2.8. (page 78) in \cite{MU4}. It presents the main properties of topological pressure function:
\begin{proposition}
The topological pressure functions satisfies the following:\\

(a)\ \ If the incidence matrix A is finitely irreducible, then for every $n \geq 1$ we have $\theta_{n,E,A} =\theta_{E,A}$\\

(b)\ \ The topological pressure function $P_{E,A}$ is non-increasing on
$[0,\infty)$, strictly decreasing on $[\theta_{E,A},\infty)$ to negative infinity,
convex and continuous on $(\theta_{E,A},\infty)$.\\

(c)\ \ $P(0)=\infty$ iff $|E|=\infty$\\

(d)\ \ $P(t)= inf \{u \geq 0 :$ 
$\exists$  $n_0 \geq 0$ so that 
$\sum_{\omega \in E_A^*;|\omega | \geq n_0} ||\varphi_{\omega}'||^t e^{-u|\omega|} < \infty \}$ 

\end{proposition}
\begin{remark}
We would like to add that if the incidence matrix is not finitely irreducible, part a) is not necessarily satisfied.\\
Consider the standard continued fractions conformal iterated function system $S_{\mathbb{N}}$, where $\mathbb{N}$ is the set of all edges (see \cite{MU3}). Define the edge incidence matrix the following way:\\
$A(i,j)=1$ iff $|i-j| \leq 1$.\\ 
The incidence matrix $A$ is irreducible, but not finitely irreducible.\\
We have: $\theta_{n,\mathbb{N},A} = \frac{1}{2n}$. Thus, the finiteness parameter for the pressure function $\theta_{\mathbb{N},A}$ is zero.
\end{remark}
\begin{definition}
We call a system $S$ \emph{regular} if there exists a $t \geq 0$ so that $P_{E,A}(t)=0$.
\end{definition}
Next, we collect the main results about finitely irreducible CGDMS's. In \cite{MU4}, these results appear in the context of a finitely primitive CGDMS.
They can be easily addapted to the case when the incidence matrix is finitely irreducible.
\begin{theorem}
A finitely irreducible CGDMS is regular if and only if there exists a t-conformal measure.
If such a measure exists, then the measure is unique and $P_{E,A}(t)=0$.
\end{theorem}
\begin{lemma}
If a finitely irreducible CGDMS is regular and $m$ is the unique h-conformal measure, then:
\[
M=min\{m(X_v):v \in V \} > 0.
\]
\end{lemma}
\begin{lemma}
If $S$ is a regular CGDMS with incidence matrix $A$ finitely irreducible, then for every $n \geq 1$ we have:
\[
1 \leq \sum_{\omega \in E_A^n} ||\varphi_{\omega}'||^h \leq K^h M^{-1}
\]
where $h$ is the unique zero of the pressure function $P_{E,A}$ and $M$ is the constant coming from
Lemma 3.1.
\end{lemma}
\begin{theorem}
Given a conformal GDMS with incidence matrix finitely irreducible, we have:
\[
HD(J_{E,A})= inf\{ t \geq 0: P_{E,A}(t) < 0 \}  = sup\{HD(J_{E,F}) | F \subset E, |F| < \infty \} \geq \theta_{E,A}(S)
\]
where $HD(J_{E,A})$ is the Hausdorff dimension of the limit set $J_{E,A}$.\\
\\
If $P_{E,A}(t)=0$, then $t$ is the only zero of the function $P_{E,A}$ and $t=HD(J_{E,A})$.
\end{theorem}
\begin{definition}
A conformal GDMS is said to be \emph{strongly regular} if there exists $t \geq 0$
such that $0 < P_{E,A}(t) < \infty$. If a conformal GDMS is not regular
we call it \emph{irregular}.
\end{definition}
An example of a regular system that is not strongly regular can be found in \cite{MU3}\\
An example of an irregular system is given in \cite{MU4}.
\newpage
\section{General Graph directed Markov Systems}
\begin{definition}
Let $S_E$ be a finite CGDMS. ($|E| < \infty$). $C \subset E$ is called a strongly connected component if for any two edges
$c_1,c_2$ from $C$ there exists $\omega \in E_A^*$ so that $c_1 \omega c_2$ is in $E_C^*$ and is the largest in the sense of inclusion.\\
$C$ is called a maximal strongly connected component if $HD(J_{C,A})=HD(J_{E,A})$. 
\end{definition}
\begin{theorem}
Given a finite CGDMS $S_E$, 
\[
HD(J_{E,A})=max\{HD(J_{C,A})| C\ \ is\ \ a\ \ strongly \ \ connected \ \ component \ \ of\ \  E\ \ \}
\]
In particular, the Hausdorff dimension of the limit set is the zero of the topological pressure.
\begin{proof}

Let $C_1,C_2,...,C_k$ be the strongly connected components and fix $1 \leq j \leq k $.\\
Any admissible word (finite or infinite) that contains at least one letter (edge) from $C_j$ can be written as
$\beta \alpha_j \gamma$,
where none of the letters of $\beta$ and $\gamma $ are from $C_j$. ($\beta $ or $\gamma $ or both may be the empty word) and 
$\alpha_j$ is the maximal string that has letters only from $C_j$.\\
If $\beta \alpha_j \gamma $ is an admissible word, then $\gamma \alpha_j \beta$ can't be an admissible word since $C_j$ is a strongly connected component.\\
Thus, we can write the following inequality:
\[
\sum_{\omega \in E_A^*} ||\varphi_{\omega}'||^t e^{-u|\omega|} \leq T \sum_{\omega_1 \in C_1^*} ||\varphi_{\omega_1}'||^t e^{-u|\omega_1|}
\sum_{\omega_2 \in C_2^*} ||\varphi_{\omega_2}'||^t e^{-u|\omega_2|}...\sum_{\omega_k \in C_k^*} ||\varphi_{\omega_k}'||^t e^{-u|\omega_k|}
\]
where the constant $T$ depends only on the set of isolated edges (edges that are not in a strongly connected component).
\\
This shows that $P_{A,E}(t)\leq max_{1 \leq j \leq k }\{P_{C_j}(t)\}$ (this follows from Proposition 3.1.).\\
Since the other inequality is obvious, we conclude that:
$P_E(t) = max_{1 \leq j \leq k }\{P_{C_j}(t)\}$.\\
Thus, using Proposition 3.1. d) and Theorem 3.2., we have:
\[
HD(J_S)=max\{HD(J_C)| C\ \ is\ \ a\ \ strongly \ \ connected \ \ component \ \ of\ \  E_A\ \ \}
\]
\end{proof}
\end{theorem}
\begin{lemma}
Let $S$ be a finite CGDMS and $E$ be the set of all edges. Let $C_1, C_2,...,C_k$ be the maximal strongly connected components.
Suppose that for every $1\leq i \ne j \leq k$, $C_i$ and $C_j$ do not communicate (there is no admissible word with letters from both components). 
Let $C_0=E \setminus \cup_{1\leq i \leq k}C_i$.\\
Then there exists $k_0 \geq 1$ and $0 < a <1$ so that $\forall n\geq 1$,\\
$Z_{C_0,n}(h)=\sum_{|\omega |=n,\omega \in C_0^*} ||\varphi_{\omega}'||^h \leq k_0 a^n$\\
In particular, there exists $M_0 > 0$ so that $\sum_{\omega \in C_0^*} ||\varphi_{\omega}'||^h \leq M_0 $ 
\begin{proof}
This is obvious since we can't have arbitrarily long words with letters only from $C_0$.
\end{proof}
\end{lemma} 
\begin{lemma}
Let $S$ be a finite CGDMS and $E$ be the set of all edges. Let $C_1, C_2,...,C_k$ be the maximal strongly connected components.
Suppose that for every $1\leq i \ne j \leq k$, $C_i$ and $C_j$ do not communicate. 
Let $C_0=E \setminus \cup_{1\leq i \leq k}C_i$.\\
For every $1\leq i\leq k$, let $C_i^{**}$ be the set of all the finite words with at least one letter form $C_i$.\\
Then there exists $M_i>0$ so that:\\
 $\sum_{|\omega |=n,\omega \in C_i^{**}} ||\varphi_{\omega}'||^h \leq M_i$ $\forall n\geq 1$.
\begin{proof}
Using Lemma 3.2. there exists $N_i>0$ so that\\
$\sum_{|\omega |=n,\omega \in C_i^*} ||\varphi_{\omega}'||^h \leq N_i$ $\forall n\geq 1$.\\
For every $0 \leq j\leq n$ let $C_{i,j}^{**}$ be the set of all the words in $C_i^{**}$ containing exactly $j$ letters from $C_i$.  
\\
Every word in $C_{i,j}^{**}$ is of one of the following forms: $\alpha \omega $,\ \ \ ($\alpha \in C_0^{n-j},\omega \in C_i^j$),
\\ $\omega \beta $,\ \ \ ($\beta \in C_0^{n-j},\omega \in C_i^j$) or \ \ \ $\alpha \omega \beta $,\ \ \ ($\alpha , \beta  \in C_0^*,\omega \in C_i^j, |\alpha |+|\beta|=n-j $).
\\Based on this observation and on Lemma 4.1., we have:\\
$\sum_{|\omega |=n,\omega \in C_{i,j}^{**}} ||\varphi_{\omega}'||^h \leq (n-j+1) k_0^2 a^{n-j} N_i $,\\
where $a$ and $k_0$ come from the Lemma 4.1.\\
Thus: $\sum_{|\omega |=n,\omega \in C_i^{**}} ||\varphi_{\omega}'||^h \leq N_i k_0^2 (1 + 2a + ... + na^{n-1})$.\\
Letting $M_i= N_i k_0^2 (1 + 2a + ... + na^{n-1} + ...)$, we get:\\
$\sum_{|\omega |=n,\omega \in C_i^{**}} ||\varphi_{\omega}'||^h \leq M_i$ $\forall n\geq 1$ 
\end{proof}
\end{lemma}
\begin{proposition}
Let $S$ be a finite CGDMS and $E$ be the set of all edges. Let $C_1, C_2,...,C_k$ the maximal strongly connected components and let $C_0=E \setminus \cup_{1\leq i \leq k}C_i$.
If for every $1\leq i \ne j \leq k$, $C_i$ and $C_j$ do not communicate, then
there exists a finite constant $M$ so that:\\
\\
 $Z_{n,E}(h)=\sum_{|\omega |=n} ||\varphi_{\omega}'||^h \leq M$, \ \ \ \ \ \ \ $\forall n\geq 1$. 
\begin{proof}
Since the maximal strongly connected components are pairwise non-communicating, we can write the following:\\
$Z_{n,E}(h)=\sum_{|\omega |=n} ||\varphi_{\omega}'||^h$ = $\sum_{|\omega |=n, \omega \in C_1^{**}} ||\varphi_{\omega}'||^h$ + ... +
$\sum_{|\omega |=n, \omega \in C_k^{**}} ||\varphi_{\omega}'||^h$ + $\sum_{|\omega |=n, \omega \in C_0^*} ||\varphi_{\omega}'||^h$ $\leq
M_1 + ... + M_k + M_0$.\\
Choosing now $M=M_1 + ... + M_k + M_0$ we get:\\
$Z_{n,E}(h)=\sum_{|\omega |=n} ||\varphi_{\omega}'||^h \leq M$, \ \ \ \ \ \ \ $\forall n\geq 1$ 
\end{proof} 
\end{proposition}
\begin{proposition}
Let $S$ be a finite CGDMS and $E$ be the set of all edges. Let $C_1, C_2,...,C_k$ the maximal strongly connected components. Let $h=HD(J_{E,A})$.
If for every $1\leq i \ne j \leq k$, $C_i$ and $C_j$ do not communicate, then $H^h(J_{E,A}) < \infty$.
\begin{proof}
For every $n\geq 1$, $\cup_{|\omega|=n}\varphi_{\omega}(X_{t(\omega)})$ is a cover of $J$ whose diameter converges to 0 as $n\rightarrow \infty$.\\
We have the following: $\sum_{|\omega |=n} diam(\varphi_{\omega}(X_{t(\omega)})^h \leq D^h \sum_{|\omega |=n}||\varphi_{\omega}'||^h < D^h M$,\\
where $M$ is a constant coming from Proposition 4.1. and $D$ is a constant coming from (4.20.), page 73 in \cite{MU4}.
Thus: $H^h(J_{E,A}) < \infty$.
\end{proof}
\end{proposition}
\begin{proposition}
Let $S$ be a finite CGDMS. Let $E$ be the set of all edges. Let $C_1, C_2,...,C_k$ the maximal strongly connected components. If there exists $1\leq i \ne j \leq k$ so that $C_i$ and $C_j$ communicate, then $sup_{n\geq 1}Z_{n,E}(h)= \infty$, where $h=HD(J_E)$.
\begin{proof}
We may assume that $C_1=\{e_{1,1},e_{1,2},...,e_{1,p}\}$ and $C_2=\{e_{2,1},e_{2,2},...,e_{2,r}\}$ are the 2 maximal strongly connected component that communicate.\\
We may also assume that there exists $\omega_0 \in [E \setminus(C_1 \cup C_2)]_A^*$ so that $e_{1,p}\omega_0 e_{2,1} \in E_A^*$.\\
For every $n \geq 1$ and $1 \leq l \leq p$, let $C_{A,l}^{1,n}$ be the set of all admissible words with letters from $C_1$, of length $n$ and whose last letter is $e_{1,l}$.\\
Define: $Z_{n,C_1,A,l}(h)= \sum_{\omega \in C_{A,l}^{1,n}} ||\varphi'_{\omega}||^h$.\\
We will prove that for every $1 \leq l \leq p$, $inf_{n \geq 1}Z_{n,C_1,A,l}(h) > 0$.\\
First we prove that there exists $1 \leq l_0 \leq p$ so that $inf_{n \geq 1}Z_{n,C_1,A,l_0}(h) > 0$.\\
For every $n \geq 1$, there exists $M_1$ (coming from Lemma 3.2.) so that:\\
$1 \leq Z_{n,C_1,A}(h) \leq M_1$.\\
Thus there exists $1 \leq l_0 \leq p$ and an increasing sequence of positive integers $i_1,i_2,...,i_j,...$ so that:\\
$inf_{j \geq 1}Z_{i_j,C_1,A,l_0}(h) > 0$\\
For every $n \geq 1$ and $j$ so that $i_j <n$ we have:\\
$Z_{n,C_1,A,l_0}(h)n\leq Z_{n-i_j,C_1,A,l_0}(h) Z_{i_j,C,A,l_0}(h) \leq M_1 Z_{i_j,C,A,l_0}(h)$.\\
This proves that $inf_{n \geq 1}Z_{n,C_1,A,l_0}(h) > 0$.\\
Since the incidence matrix $A$ restricted to $C_1$ is finitely irreducible, we have:
\[
\mu_1 = inf_{n \geq 1, 1 \leq l \leq p}Z_{n,C_1,A,l}(h) > 0
\] 
Similarly, for every $n \geq 1$ and $1 \leq m \leq r$ we define $C_{m,A}^{2,n}$ to be the set of all the admissible words of length $n$ with letter only from $C_2$ and whose first letter is $e_{2,m}$.\\
Let $Z_{n,C_2,A,m}(h)= \sum_{\omega \in C_{m,A}^{2,n}} ||\varphi'_{\omega}||^h$.\\
We also have:
\[
\mu_2 = inf_{n \geq 1, 1 \leq m \leq r}Z_{n,C_2,A,m}(h) > 0
\]
For every finite admissible word of the form $\alpha e_{1,p}\omega_0 e_{2,1} \beta$ we have:
\[
||\varphi_{\alpha e_{1,p}\omega_0 e_{2,1} \beta}'|| \geq K^{-2} ||\varphi_{\omega_0}'||\ \ ||\varphi_{\alpha e_{1,p}}'|| \ \ 
||\varphi_{e_{2,1} \beta}'||
\]
Thus we can write the following inequality :
\[
Z_{n+|\omega_0 |,E,A}(h) \geq K^{-2h} ||\varphi_{\omega_0}'||^h \mu_1 \mu_2 n
\]
and this finishes the proof.

\end{proof}
\end{proposition}
\begin{proposition}
 Let $S$ be a finite CGDMS. Let $E$ be the set of all edges. Let $C_1, C_2,...,C_l$ the maximal strongly connected components. If there exists $1\leq i \ne j \leq l$ so that $C_i$ and $C_j$ communicate, then $H^h(J_{E,A})= \infty$, where $h=HD(J_{E,A})$.
\begin{proof}
Without loss of generality we may assume that $E= C \cup F \cup G$, where $C= \{c_1,c_2,...,c_k \}$ and $F=\{f_1,f_2,....f_p\}$ are two maximal strongly connected components, $G=\{g_1,g_2,...,g_r\}$ is a set of isolated edges and $A(c_k,f_1)=A(f_p,g_1)=1$.\\
Also, $A(g_i,g_{i+1})=1$ for every $1 \leq i \leq r-1$. 
For every $b \in F$, let $J_{b,F} = \{\pi(b\omega)|  \omega \in E_F^{\infty} \}$.\\
We have $J_F=\cup_{b \in F} J_{b,F}$ and 
$H^h(J_F) > 0$.\\
Thus there exists $b_0 \in F$ so that $H^h(J_{b_0,F}) > 0$.\\
Let $b_1 \ne b_0$ in $F$ so that $A(b_1,b_0)=1$. Let $g=(g_1,g_2,...,g_r)$.\\
For $n \geq 0$ large enough, let:\\
$T_n=\{\omega = \beta g \gamma b_1 \ \ | \ \ |\omega|=n, \beta \in E_C^*,\gamma \in E_F^*\}$.
We have: $J_{E,A} \supset \cup_{\omega \in T} \varphi_{\omega}(J_{b_0,F})$.\\
For every $a,b \in C$, $H^h(\varphi_a(X_{t(a)})\cap \varphi_b(X_{t(b)})=0$. For every $a_1,b_1 \in F$, $H^h(\varphi_{a_1}(X_{t(a_1)})\cap \varphi_{b_1}(X_{t(b_1)})=0$.\\
Thus, $H^h(\varphi_{\omega}(X_{t(\omega)}) 
\cap \varphi_{\rho (X_{t(\rho)}})=0$, 
for every $\omega, \rho \in T_n$.\\
So, $H^h(J_{E,A}) = \sum_{\omega \in T_n} H^h(\varphi_{\omega}(J_{b_0,F})) 
\geq \sum_{\varphi_{\omega} \in T_n} 
K^{-1} || \varphi_{\omega}'||^h 
H^h(J_{b_0,F})$.\\
As in the previous theorem, $sup_n\{\sum_{\varphi_{\omega} \in T_n} ||\varphi_{\omega}'||^h \} =\infty$.\\
In conclusion, $H^h(J_{E,A}) = \infty$.    
\end{proof}
\end{proposition}
\begin{theorem}
Let $S$ be an infinite CGDMS and let $E=\mathbf{N}$ be the set of all edges. Suppose that for any two vertices $v_1,v_2$ the exists an edge from $v_1$ to $v_2$. If the incidence matrix $A_E$ is irreducible, then:
\[
HD(J_S)= sup \{HD(J_F)| F \subset E \ \ finite \} = inf \{ t|P_{E,A}(t)<0 \}
\]
\begin{proof}
Since the incidence matrix $A_E$ is irreducible, there exists an increasing sequence of finite sets of edges
$\{E_k\}_{k \geq 1}$ so that $\cup_{k \geq 1 }E_k = E$ and $A_{| E_k \times E_k}$ is irreducible. Relabeling the edges
by the positive integers, we can assume that there exists an increasing sequence of positive integers $N_1,N_2,...,N_i,...$ so that for every 
$i$, $A|_{\{1,2,...,N_i\}^2}$ is irreducible.\\
For every $i \geq 1$, let $F_i= \{1,2,...,N_i\}$. Using Lemma 3.2., for every $i \geq 1$ there exists a constant $M_i>0$ so that for every $n \geq 1$:\\
\[
1 \leq 
\sum_{\omega \in F_i^n} ||\varphi '||^{h_i} \leq K^{h_i} M_i
\]
where $h_i$ is the Hausdorff dimension of $J_{F_i}$.\\
In particular, for every $i$, $M_i = min_v\{m_i(X_v)\} > 0 $, where $m_i$ is the conformal measure on $J_{F_i}$ (Lemma 3.1.).\\
Next, we will show that $inf_i\{M_i\}>0$. \\
Let $v_1,v_2,...,v_r$ be the set of vertices. For any two vertices $v_p,v_s$ there exists a finite path $\omega_{ps}$ so that $\varphi_{\omega_{ps}}(X_{v_p}) \subset X_{v_s}$.\\
Choose $i_0$ large enough so that for every $1 \leq p,s \leq r$, $\omega_{ps} \in F_{i_0}^*$.\\
Let $\epsilon > 0$ arbitrarily small. Suppose there exists $i>0$ and $s$ so that $m_i(X_{v_s}) < \epsilon$.\\
Then for every other vertex $v_p$ we have:\\
$\varphi_{\omega_{ps}}(X_{v_p}) \subset X_{v_s}.$ Thus $m_i(\varphi_{\omega_{ps}}(X_{v_p})) < \epsilon$.\\
But $m_i(\varphi_{\omega_{ps}}(X_{v_p})) = \int_{X_{v_p}}{|\varphi_{\omega_{ps}}'|}^{h_i} dm_i$.\\
So $m_i(X_{v_p}) \leq K^{h_i} ||\varphi_{\omega_{ps}}'||^{-h_i} \epsilon. \leq K^d z^d \epsilon$,  \\
where $z= inf_{p,s}|| \varphi_{\omega_{ps}}' ||$.
This is a contradiction though, since $\sum_p m_i(X_{v_p})=1$.\\
Thus there exists a constant $M$ so that for every $n \geq 1$ and $i \geq 1$:
\[
1 \leq \sum_{\omega \in E_i^n} ||\varphi_{\omega}'||^{h_i} \leq K^{h_i} M
\]
Now we will follow the ideas from the proof of Theorem 4.2.13 in \cite{MU4}.\\
Let $\eta = sup \{HD(J_F)| f \subset E finite\}$. Let $t> \eta$.\\
For every $n \geq 1$:\\
$\sum_{\omega \in E^n} ||\varphi_{\omega}'||^t = sup_i \sum_{\omega \in E_i^n} ||\varphi_{\omega}'||^t \leq
sup_i \sum_{\omega \in E_i^n} 
|| \varphi_{\omega}'||^{h_i} 
s^{n(t-h_i)}$\\ 
$\leq
s^{n(t-\eta)} 
\sum_{\omega \in E_i^n} ||\varphi_{\omega}'||^{h_i} \leq s^{n(t-\eta)} K^{\eta} M.$\\
Thus $P_E(t) \leq (t - \eta) \ \ \ln s < 0$ and this finishes the proof. 
\end{proof}
\end{theorem}  
\section{Examples}
We would like to end with an example showing that Theorem 3.2. doesn't necessarily hold for a CGDMS whose incidence matrix is not assumed to be finitely irreducible.\\
\\
Consider the standard continued fraction CIFS $\{S_{\mathbb{N}}\}_{e \in \mathbb{N}}$. There is only one vertex $v$, and $X_v = [0,1]$.
\\
The set of edges is the set of positive integers, and for every $e \geq 1$, $\varphi_e(x) = \frac{1}{e+x}$\\
\\
Now let us define the incidence matrix.\\
We have $A(e,f)=1$ iff $e<f$.\\
Any finite subset of the positive integers generates an empty limit set since we can't have arbitrarily long admissible words.\\
The finiteness parameter for our system is $\frac{1}{2}$ and thus $\frac{1}{2} \leq inf\{t| \ \ P_{\mathbb{N},A} < 0\}$.\\
In the same time, $sup \{HD(J_F) \ \ | \ \ F finite \}=0$.\\
This proves that Theorem 3.2 doesn't hold in this case.

\end{document}